\documentclass[11pt,oneside]{article}
\usepackage[margin=1in]{geometry}
\usepackage{verbatim,wrapfig}
\usepackage{graphicx,caption}
\usepackage{amssymb}
\usepackage{amsthm}
\usepackage{latexsym}
\usepackage{amsmath}

\usepackage{amsfonts}
\usepackage[colorlinks=true]{hyperref}

\numberwithin{equation}{section}
\numberwithin{figure}{section}

\makeatletter
\g@addto@macro\@floatboxreset\centering
\let\@fnsymbol\@alph
\makeatother

\newcommand{\fulline}{ \begin{wrapfigure}{R}{0\textwidth}\centering\end{wrapfigure}} 

\title{Inquiry-based learning in a first-year honors course}

\begin{document}

\author{Diana Davis
}
\maketitle

\begin{abstract}
We describe a case study of a  problem-solving section, using the ``Harkness'' discussion method, of an honors multivariable calculus course. Students in the problem-solving section had equivalent outcomes  on exams,  reported higher ratings in self-assessments of skills, and took more math classes in the following year, compared to students in the lecture-based sections.
\end{abstract}


\section{Introduction}

When talking about ``alternative'' methods of teaching, we would do well to remember that lecturing is only standard because of a historical anachronism. Lecturing started before the printing press was invented, before universal literacy, when the one person who could read stood in front of those who could not, and read aloud the one available copy of the book. 

In 2017, any student can learn college-level mathematics for free from online texts, or from online lectures such as Khan Academy or MIT OpenCourseWare. If we want students to feel that there is a reason to pay for and attend college classes, there has to be an added value from the student's actual presence in the classroom. A major benefit of students and the instructor all being in the same room is that they can learn by talking with each other, and this idea has gained wide popularity, with a nationwide push towards ``active learning'' in the college math classroom.

Active learning appeared in mathematics courses more than a century ago, when R.L. Moore developed his ``Moore method'' between 1911-1916 \cite{moore}. Under the Moore method, students develop and present key mathematical ideas before their peers, and the instructor acts as a mentor and guide, aiding in the students' own discovery. In order for this and other ``inquiry-based'' methods to work optimally, students need to ask each other good questions, which they often have never learned to do, so some instructors have tried intentionally modeling asking good questions as scaffolding to guide students to ask them on their own \cite{belcastro}.

Many instructors have found success using the ``flipped classroom'' model, where students prepare for class by learning the material on their own via textbook, slides or video lectures, and then spent class time solving math problems, playing content-specific games, or doing computer programming \cite{love}.
For a comparison of the ``flipped classroom'' versus ``inquiry-based learning'' teaching methods, and methods of combining the two, see \cite{capaldi}. 
Instructors have found considerable success using class time for student problem-solving in small groups, both in student pass rates and in student satisfaction; see \cite{deshler}.
There are also many ways to incorporate small amounts of active learning into an existing lecture class; see \cite{marotta}.

Here we discuss an inquiry-based learning approach called the \textit{Harkness method}; for its history and a discussion of its use in classes in many disciplines, see \cite{book}. We also use a problem-based curriculum approach that interweaves many different topics and ideas into a single assignment; for the history of this approach and a description of an example class, see \cite{harkness}.

Some students don't like inquiry-based learning. Well, I never liked lectures! After a wonderful mathematical experience in high school exploring problems and discussing them every day using this method, I headed to a top liberal arts college, six of whose mathematics faculty members have won the highest national award for their teaching
$-$ where I sat unhappily through four years of doubtless exceptional lectures. Many students and faculty members would prefer to teach and learn in an inquiry-based method, rather than being forced to engage in the lecture method, and they should be encouraged to do so. This paper describes a case study of the \textit{Harkness method} approach to inquiry-based learning, and gives evidence that it was successful.

\section{Structure of the experiment}\label{structure}
In spring 2016 in the math department at Northwestern University, we ran an experimental section of the third quarter of our first-year honors sequence, Math 290. All of the instructors for the honors sequence were interested in using active learning methods, and the department agreed to allow a one-quarter ``pilot'' inquiry-based section. The students in Math 290 were a relatively homogeneous population, because Northwestern sections students into one of four different linear algebra and multivariable calculus tracks. The Math 290 track was for students who scored a 4 or 5 on the AP calculus exam, but who were not ready for an intensive proof-based course. Additionally, there was some attrition over the course of the year, so the students who took Math 290-3 in the spring were the best of the Math 290 group.

In the winter quarter, all students were given information about this ``problem-solving section,'' and had the opportunity to sign up: Every student who wanted to be in the problem-solving section was allowed to register for it, and any student who didn't want to be in it registered for one of the two lecture-based sections. Thus registration for the experimental section was on an opt-in basis, in order to give the course the best possible chance of success and also to be fair to the students. We had $10$ students in the problem-solving section, and $10$ and $13$ students in the two lecture-based sections.

We measured the exam performance of students in all three sections, and also surveyed all of the students at the beginning, middle and end of the course.\footnote{The teaching method was within the scope of standard pedagogy, and the questionnaires were within the scope of standard course evaluations, so we did not need to seek IRB approval for human subjects.}  There was no statistically significant difference in exam scores under the two teaching methods (see $\S\ref{exams}$). On the surveys, students reported how much they liked the course, how comfortable they were discussing math with others, and how comfortable they were speaking in front of others. On each of these measures, the self ratings by students in the problem-solving section increased more, and ended up  higher (better), than those of students in the lecture-based sections (see $\S\ref{surveys}$). We also asked several other questions to gauge the students' experience with the different teaching methods. Finally, a year after the conclusion of the course, we gathered data from the Registrar on how many math courses each student had taken in the subsequent (sophomore) year, and how many of the students had declared the math major (see $\S$\ref{longrange}).

\newpage
\section{Instructors}\label{instructors}

All of the instructors for this course in the fall, winter and spring were Postdoctoral Lecturers, faculty members on a three-year contract who were hired based primarily on teaching ability. 

 As we compare outcomes of students under the problem-based verses lecture-based methods, it is important to note that the lecturers were excellent. In student evaluations at Northwestern, instructors are graded on a scale from $1$ to $6$ on six different measures (preparedness, effective communication, enthusiasm, etc.), and each of the two lecturers averages over $5.00$ overall, which are among the highest in the Northwestern math department. It should also be emphasized that {\it all} instructors included some active learning in their ``lecture,'' for example asking for student input or having students work out examples together during class.

\section{Design of the problem-based curriculum}

The problems that formed the content of the course were based on the math curriculum at Phillips Exeter Academy. About half of the problems were taken directly from the PEA {\it Mathematics 5} book \cite{math5}, some of the problems were taken from  existing Math 290 lectures and problem sets, and the rest I wrote specifically for the course. All of the content of the course was introduced through the problems: there were no lectures, no worked examples, no summarizing content except what was worked into the text of each problem. The problems were designed to introduce the material $-$ definitions, theorems and examples $-$ by starting with simple problems and building in difficulty and complexity, with a mix of different ideas every day. For more information on the Exeter math curriculum design, see \cite{harkness}. 

To construct the curriculum, I made a list of the topics covered in the course and considered the type of problems on each topic that the students would need to be able to solve on the midterms and final exam. Then I wrote a sequence of problems that would gradually introduce the ideas, starting with problems that introduced basic definitions and culminating in the type of problems they would need to do on exams. I wrote or typed each problem on a strip of paper (see picture below left). Then I laid out a blank piece of paper for each class day, and arranged the daily assignments (see picture below right). Given the sequence of problems building the ideas of a particular topic, I would usually put one problem on that topic per day, with a gap of two days if the idea needed more time to be absorbed, and two problems building the same topic on a single homework assignment if they used different ideas. The full text is at \cite{text}.

\begin{figure}[!h]
\noindent \includegraphics[trim={0 2.7cm  0 0},clip,width=0.455\textwidth]{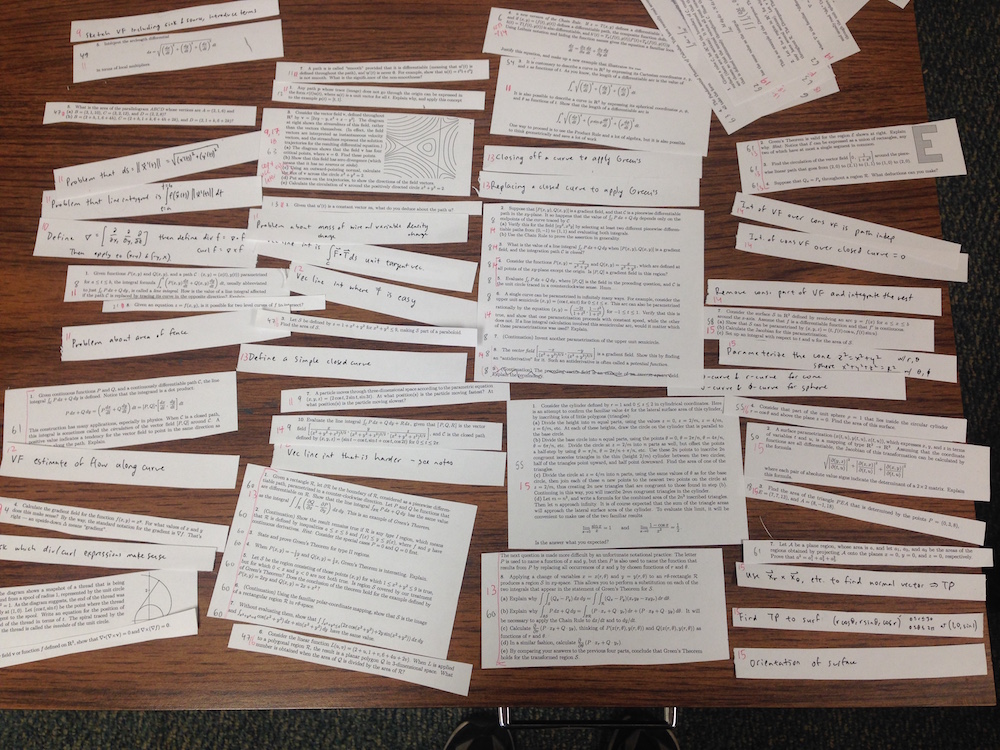} 
\includegraphics[trim={0 1cm  0 5cm},clip,width=0.53\textwidth]{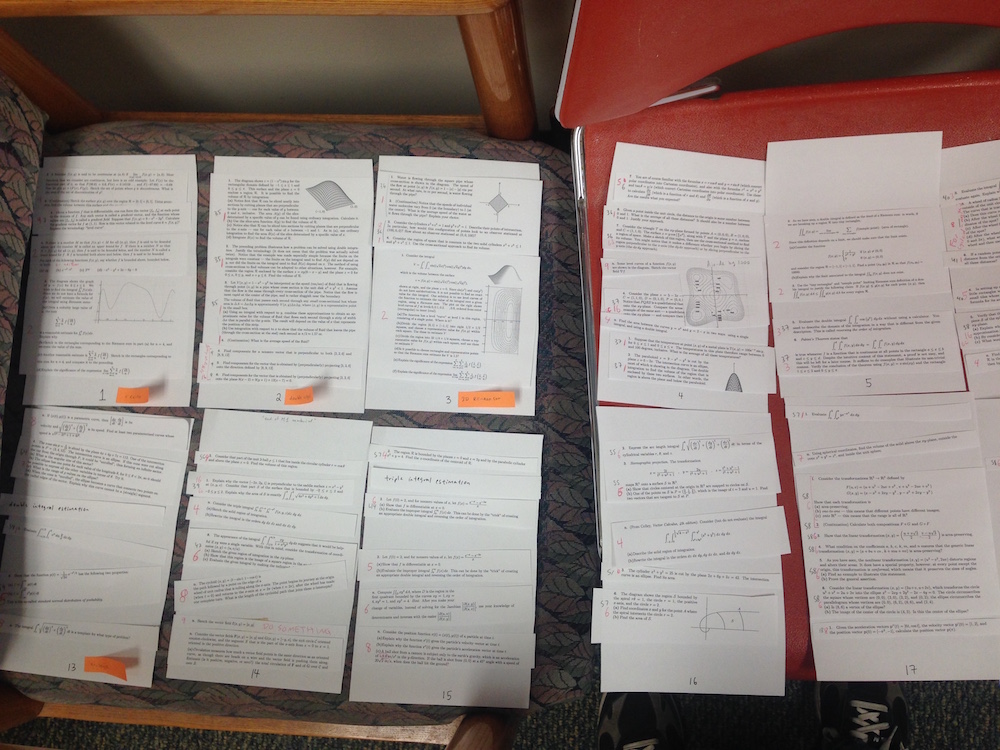} 
\captionsetup{labelformat=empty}
\caption{The problems, arranged by topic (left) and by daily assignment (right)}
\end{figure}

The following example shows the sequence of problems on Green's Theorem. 

\vspace{0.5em} \noindent {\bf Green's Theorem:} {\it If $D$ is a closed, bounded region in $\mathbf{R}^2$ whose boundary $C$ consists of finitely many simple, closed, piecewise-differentiable curves, oriented so that $D$ is on the left when one traverses $C$, and if $\mathbf{F} = [P,Q]$ is differentiable everywhere in $D$, then }
\[ \oint_C P \ dx + Q \ dy = \iint_D \left(\frac{\partial Q}{\partial x} - \frac{\partial P}{\partial y}\right) \ dx \ dy.
\]
{\bf On the exam, students will need to:} See a problem  about a vector line integral over a path, know to add a curve to close off the path, apply Green's Theorem to the resulting closed region and its simple closed curve boundary, and subtract off the vector line integral of the added part.

\vspace{0.5em}\noindent {\bf Sequence of problems:} \vspace{-1em}

\begin{itemize}
\item Learn to set up and compute double integrals (done in the first part of the course). \vspace{-0.5em}
\item Given a graph of level curves for a function $f(x,y)$, sketch the gradient field $\nabla f$. \vspace{-0.5em}
\item Sketch simple vector fields like $[x,y]$ and $[-y,x]$.\vspace{-0.5em}
\item Given a picture of a vector field with several directed curves, estimate whether the vector line integral over each curve is positive, negative or $0$. \vspace{-0.5em}
\item Set up and compute a given vector line integral.\vspace{-0.5em}
\item Say whether each of several given curves is simple, and whether it is closed.\vspace{-0.5em}
\item Given several directed curves, for each one shade the region ``to the left.'' \vspace{-0.5em}
\item Compute both sides of the Green's Theorem equation for a given vector field and region.\vspace{-0.5em}
\item An impossible vector line integral that is easy to compute after applying Green's Theorem.\vspace{-0.5em}
\item An impossible area integral that is easy to compute after applying Green's Theorem.\vspace{-0.5em}
\item A vector line integral over a non-closed curve that gives a hint to close off the region.\vspace{-0.5em}
\item A vector line integral over a non-closed curve with no hint.
\end{itemize}

\section{Design of the problem-based course}\label{design}

The problem-solving class met four times per week for $50$ minutes, and the format was the same every day: the students were assigned about six problems for homework, on which they were expected to spend at least $90$ minutes. 

\begin{wrapfigure}[15]{R}{0.7\textwidth}
\centering
\vspace{-1.5em}
\includegraphics[trim={8cm 0 12cm 5cm},clip,width=0.7\textwidth]{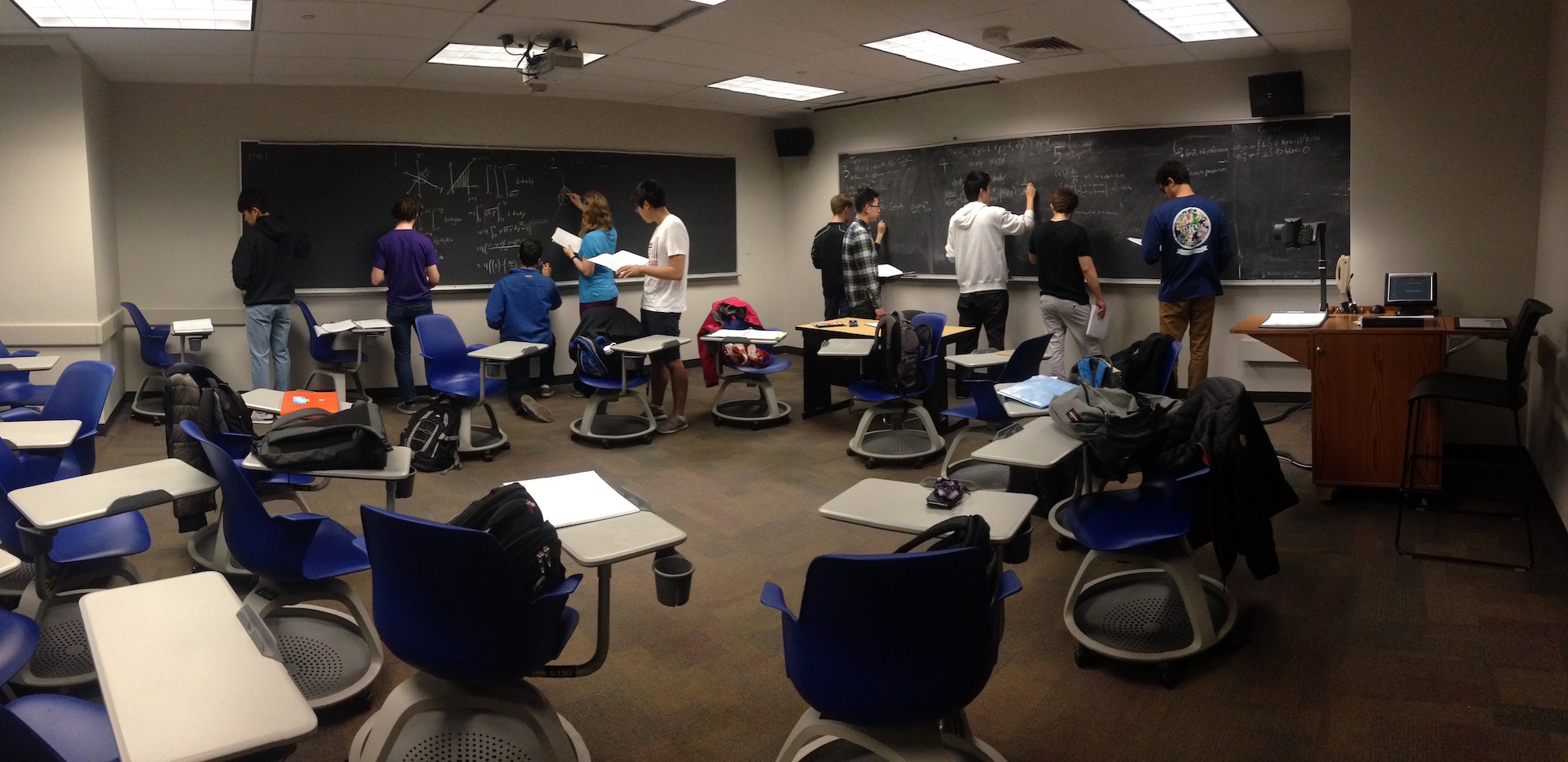}
\captionsetup{labelformat=empty}
\caption{Students writing their solutions on the board at the start of class}
\end{wrapfigure}
Before the class started, I arranged $11$ chairs in a circle, moved the teacher table and podium to the side, and wrote the problem numbers from the previous night's homework (e.g., 1-6) on the board.  Each student chose a problem to put up on a first-come, first-served basis and copied their solution from their notebook onto the board. The picture below is from the first few minutes of class, when all $10$ students have put their backpacks down on their chairs and are writing their solutions on the board, some singly and some in pairs.\footnote{A video of students writing up and discussing solutions: \url{https://www.youtube.com/watch?v=ZNFAmob_Odk} } This took place for about ten minutes, starting about five minutes before class officially started. 


By five minutes into class, students were generally clustered in pairs or groups of three, comparing their work to what was on the board and asking questions of the student who had written up the work. At this point, I would ask everyone to sit down $-$ I felt bad cutting off the excellent mathematical discussions that they were having at the board, but I wanted them to discuss those issues and have their insights with the whole class, so that everyone could benefit from the intellectual journey and not just hear the clarified answers.

We spent the remaining $45$ minutes of class sitting in a circle  discussing the problems, using the ``Harkness method'' (developed in the math department at Phillips Exeter Academy, where I have been a student and a faculty member). You can see the circle of chairs in the picture to the right. The students all opened their problem books to the previous night's problems and opened their notebooks to where they had done their work, and then the student(s) who had put up problem $\#1$ went up to the board, stood next to their writing and explained their solution. Sometimes everyone had done the same thing and there was no discussion; more frequently, discussion ensued. Students asked questions if they thought the student's solution was unclear or if they thought it was incorrect; they asked about how this problem related to previous problems; they reported when they had a different final answer and worked to resolve the discrepancy; they explained their different solution (sometimes putting their notebook on the document camera to show their work) and asked if their method was also correct. Often the discussion of a single problem would involve five or six of the ten students.\footnote{A video of students presenting and discussing a homework problem: \url{https://www.youtube.com/watch?v=BS_zwKX_6U8&t=200s} } When I decided we had discussed that problem enough, the student would sit down and the student who put up $\#2$ would go up to the board, and we would do the same thing for the next problem, and so on for all of them. 

\begin{wrapfigure}[14]{L}{0.6\textwidth}
\centering
\vspace{-1.5em}
\includegraphics[width=0.6\textwidth]{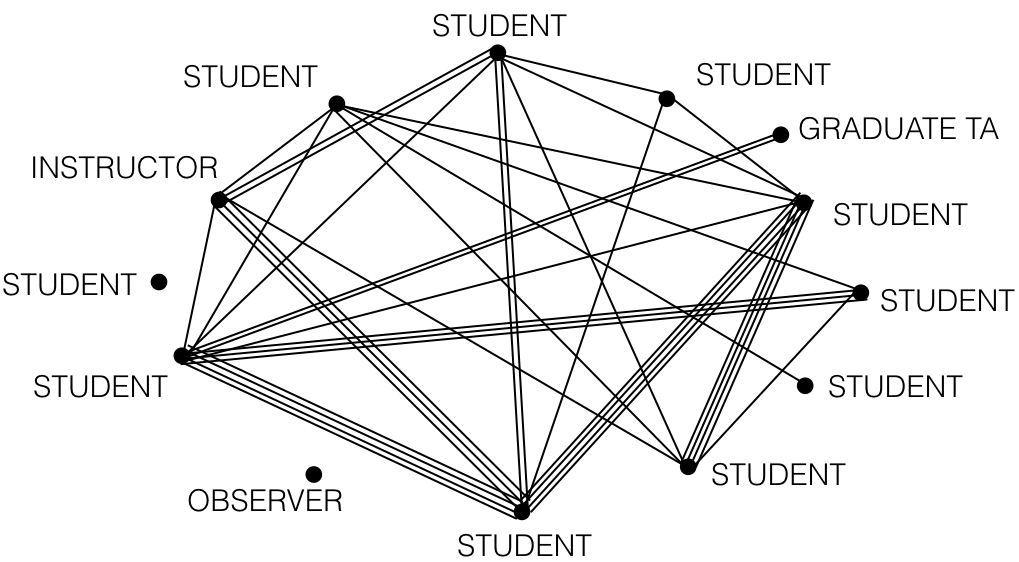}
\captionsetup{labelformat=empty}
\caption{The class discussion during one class period}
\end{wrapfigure}
The interaction diagram to the left shows the class discussion during a typical class period a few weeks into the quarter. If one person spoke, and then another person spoke, there is a line from one to the other. Multiple lines indicate multiple interactions between the same pair. People are arranged in an oval, which reflects their position around the oval of chairs in the classroom. 
The conversation is reasonably balanced, with interactions between many different pairs of students. There is no one student or instructor who is dominating the discussion. You can see that one student did not speak at all on this particular day. When I showed the students this diagram during the following class, he said, ``I am just a dot!'' and resolved to ``be a line'' in the future.

In emails before the term started, on the first day of the course, on the first pages of the problem book \cite{text}, and frequently during and at the end of classes, I emphasized the various goals of class time: to discuss solutions, ask each other questions, answer each other's questions, to discover mathematics, and for the interaction diagram of our discussion to ``look like a web.''\footnote{By comparison, a diagram of a class in the Socratic method ``looks like a fan.''}  When the students had finished discussing a problem and I felt that it had gone particularly well, I told them so, pointing out specific things that they had done well individually and as a class. The students were invested in the class and engaged in the method, and I worked to consistently reinforce their growing discussion skills and class rapport.

\section{Structure of the lecture-based comparison sections}\label{salons}

The lecture-based sections, for all students in the fall and winter quarters and for the two lecture-based sections in the spring, met Monday, Wednesday and Friday for 50 minutes, plus a 50-minute ``Salon'' on Thursdays.\footnote{The term \emph{Salon} was chosen in the sense of the  literary and philosophical gatherings in $18^\text{th}$-century  France.}

The once-weekly Salon consisted of a packet of three problems, which were chosen to be interesting and somehow beyond routine problems. Students were randomly assigned to one of three groups, and each group worked on a different problem. After about $20$ minutes, all groups were expected to have solved their problems, and the students were shuffled to new groups so that each group had one ``expert'' on each problem, and the students explained solutions to each other. (The random assignments were done by having each student draw a playing card when they entered the classroom; the number determined their first group and the suit their second.) For the last $10$ minutes, students wrote up their own solutions to each of the three problems, which they handed in for a grade. This Salon took the place of a weekly TA session and quiz, which was the standard use of this fourth class period at Northwestern.\footnote{The math department has since made a change, and now the fourth (discussion section) class period is used for active learning problem sessions in all lower-level courses!} 

So that they could also have instructor feedback and a grade on work that was handed in, students in the problem-based section in the spring had two quizzes, one halfway between the first day of class and the first midterm, and the other halfway between the two midterms.



\begin{wrapfigure}[14]{R}{0.63\textwidth}
\centering
\vspace{-4em}
\includegraphics[width=0.63\textwidth]{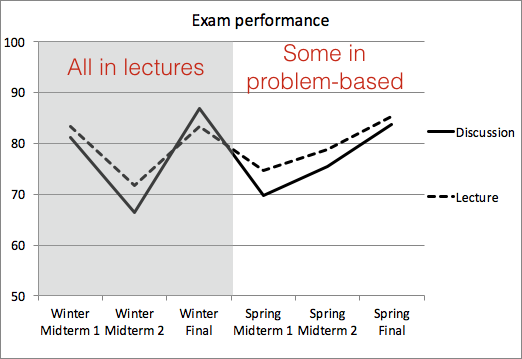}
\end{wrapfigure}
\section{Exam performance}\label{exams}
All of the students who took Math 290-3 also took Math 290-2 in the winter (and Math 290-1 in the fall) in lecture-based sections, so we were able to compare students' performance on exams under the two instruction methods. We can use the winter exam performance of students who took lecture-based sections in the spring, and of students who took the problem-solving section in the spring, as a benchmark for comparison. 

There was no systematic difference in exam scores under the two methods. The group of $10$ students in the problem-solving section averaged $1.3$ points lower than the other students across the exams in the winter (when they were in lectures) and $2.2$ points lower in the spring (when they were solving problems). This difference is not statistically significant. See table below.  

\begin{table}[!h]
\centering
\label{scores}
\begin{tabular}{c|ccc|c|ccc|c|}
 & Winter M1 & Winter M2 & Final  & Avg. & Spring M1 & Spring M2 & Final & Avg. \\ \hline
Discussion  & 81.2 & 66.4 & 86.9 &78.2&  69.8& 75.5&  83.9&76.4  \\  
Lecture  &  83.4&   71.7 & 83.3& 79.5& 74.7& 78.7&85.3 &  79.6\\ \hline
Difference & -2.2 & -5.3 & +3.6& -1.3& -4.9 &-3.2 &-1.6& -2.2
\end{tabular}
\end{table}

The data also shows that, as the term progressed, the score gap between the problem-based and lecture-based classes decreased. I believe that this is due to the ``start-up cost'' of an integrated curriculum. In the problem-based course, all of the topics studied over the term were slowly introduced starting at the very beginning, with ideas and skills slowly built up over several weeks. This meant that the students in the problem-based course spent comparatively less time studying the first midterm's material before the first midterm exam than did the students in the lecture-based sections, because they were also studying ideas that would come up later. As the term progressed, the ideas that had been slowly introduced in the problem-based course at the beginning were studied in depth and then tested in the second midterm and final exam, which may explain the narrowing of the gap. I believe that if this integrated problem-based curriculum were used consistently, any score difference would disappear, and eventually the problem-based, integrated curriculum would yield higher scores due to better retention of understanding.


\section{Survey results}\label{surveys}
All students in Math 290-3 took three surveys: on the first day of class, at the midpoint of the quarter, and in the last week. They answered the same set of five questions on all three surveys, plus some additional questions on the last survey. Every student completed all three surveys. For the first three questions, they rated how much they agreed with the given statement on a scale from $1$ to $5$. The solid line represents the average for the problem-solving section, and the dashed line represents the average across both lecture-based sections.

The $y$-axis labels are from $3.5$ to $5.0$ on all graphs, so that the trends are more visible.

\begin{wrapfigure}[12]{R}{0.5\textwidth}
\centering
\vspace{-1em}
\includegraphics[width=0.49\textwidth]{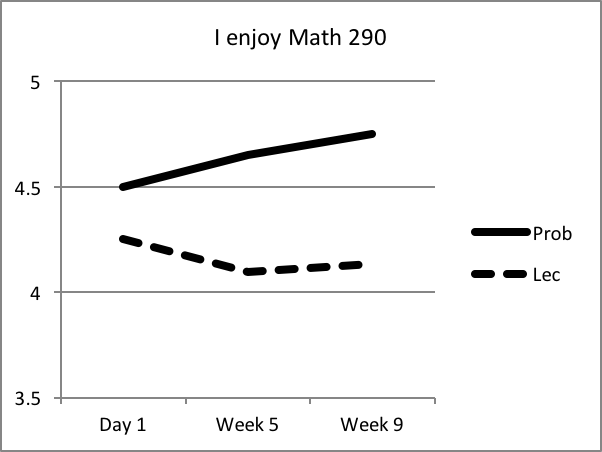}
\end{wrapfigure}
\subsection*{I enjoy Math 290} 
We asked this question because the math department was concerned about running an experiment in a course that already had very high ``customer satisfaction'': historically in this course, student ratings of instructors are high, other departments and College advisors like the course, and parents do not complain. We found that student satisfaction was high across all three sections of the course. Students enjoyed the course more (difference {\bf is} statistically significant) in the problem-solving section than in the lecture-based sections, and their enjoyment increased more (difference in increase is not statistically significant), even though all students {\it chose} whether to take a problem-solving or lecture-based course. 
\fulline

\subsection*{I feel comfortable discussing math and sharing my ideas with other students} 
\begin{wrapfigure}[14]{R}{0.5\textwidth}
\centering
\vspace{-1em}
\includegraphics[width=0.5\textwidth]{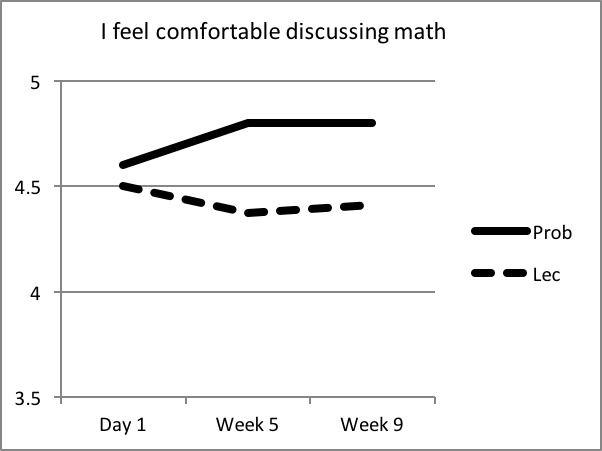}
\end{wrapfigure}
As mentioned in $\S\ref{salons}$, all students in the course had weekly ``Salons'' in the fall and winter quarters, where they worked together with other students to solve problems, and explained their solutions to others. The students in the lecture-based sections continued having weekly Salons in the spring, and students in the problem-solving section discussed problems for the entire class time four days a week, so students in all sections were discussing problems with others in the course. Unsurprisingly, then, students in all sections rated this question highly. Students in the problem-solving section reported higher self ratings, and their responses show a larger increase in average score, than those in the lecture-based sections (differences are not statistically significant).

\subsection*{I feel comfortable standing up and presenting to a group of people}
\begin{wrapfigure}[14]{R}{0.5\textwidth}
\centering
\vspace{-1em}
\includegraphics[width=0.5\textwidth]{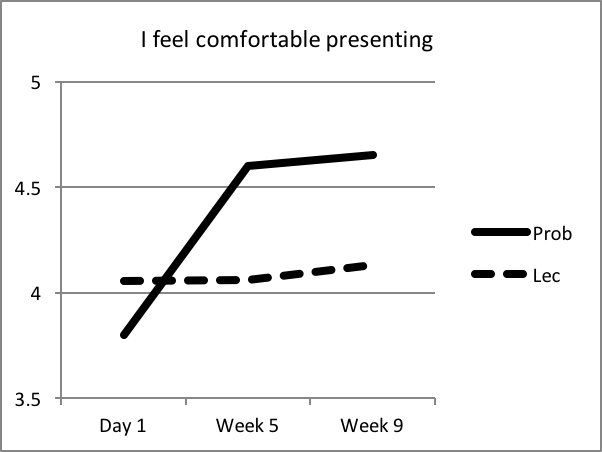}
\end{wrapfigure}
We asked this question because this is a skill that students in the problem-solving section practiced every day in class, one that would carry over into other aspects of their lives. As expected, students in the problem-solving section gave higher self ratings at the end of the course (difference is not statistically significant) and their average self ratings showed a larger increase  (difference in increase \textbf{is} statistically significant) over the two months of the course, than in the lecture-based sections, where the students' average self ratings were essentially constant.

A large proportion of Northwestern students plan to go into business, finance or consulting after they graduate, so standing up and presenting in front of other people is a highly valued skill in this population. Students in the problem-solving section improved more, and ended up better, at presenting to a group, at least as measured by their self perceptions, than students in the lecture-based section.

\clearpage
\subsection*{Given a choice between working on a problem I already know how to solve, or a problem that I don't know how to solve yet, I would rather work on$\ldots$} 
\begin{wrapfigure}[14]{R}{0.5\textwidth}
\centering
\vspace{-1.5em}
\includegraphics[width=0.5\textwidth]{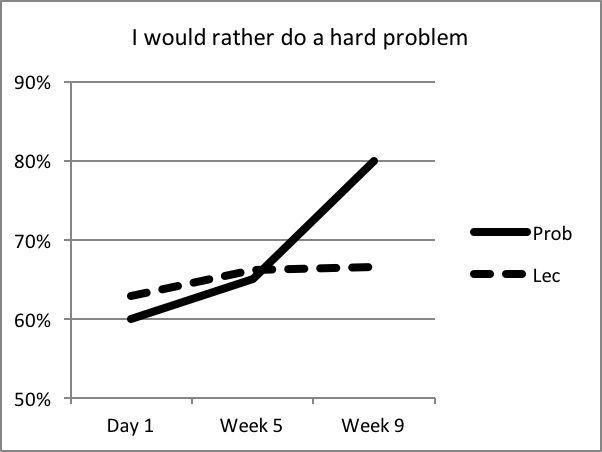}
\end{wrapfigure}
We asked this question to measure inquisitiveness and ``growth mindset.'' Students at Northwestern tend to care very much about their grades and about getting the right answer. For example, it regularly happens that a student will do all of the problems in the textbook and check their answers until they have them all correct, and then ask for more practice problems. 

We hoped that by forcing the students in the problem-solving section to struggle with hard problems, they would develop a zest for working on things that they did not know how to do. Happily, in every section of the course, more than half of the students said they would choose to work on a problem that they didn't know how to solve. The problem-solving section had a greater increase over the course of the nine weeks in the proportion of students choosing to work on the novel problem, and ended up with a higher proportion  (differences are not statistically significant).

\subsection*{During the nine weeks of the spring quarter, about how often did you miss (i.e. skip class or not attend class for any reason) your Math 290-3 class?}
\begin{wrapfigure}[11]{R}{0.6\textwidth}
\centering
\vspace{-1.5em}
\includegraphics[width=0.6\textwidth]{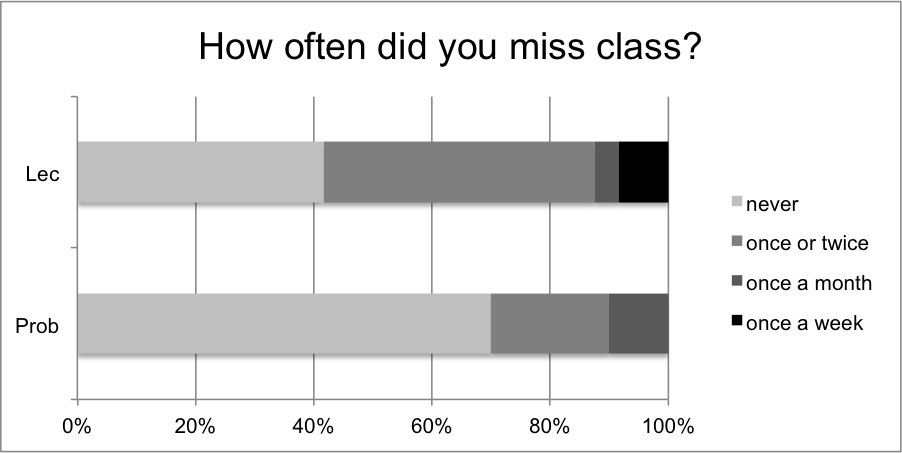}
\end{wrapfigure}
We asked this question at the end of the quarter to measure engagement in the course. Students in the problem-solving section attended class more often than students in the other sections (difference is not statistically significant). 

As the instructor, I can report that almost every day, all $10$ students were present. They also arrived early so that they could have the first choice of the problem for which they would write up the solution. In contrast, students were chronically a few minutes late for the lecture-based sections, with perhaps a third of the students arriving after the instructor had started the class.

\clearpage
\subsection*{During the spring quarter, about how often did you use an electronic device for nonacademic purposes (texting, games, etc.) during Math 290-3 class?}
\begin{wrapfigure}[11]{R}{0.6\textwidth}
\centering
\vspace{-1.5em}
\includegraphics[width=0.6\textwidth]{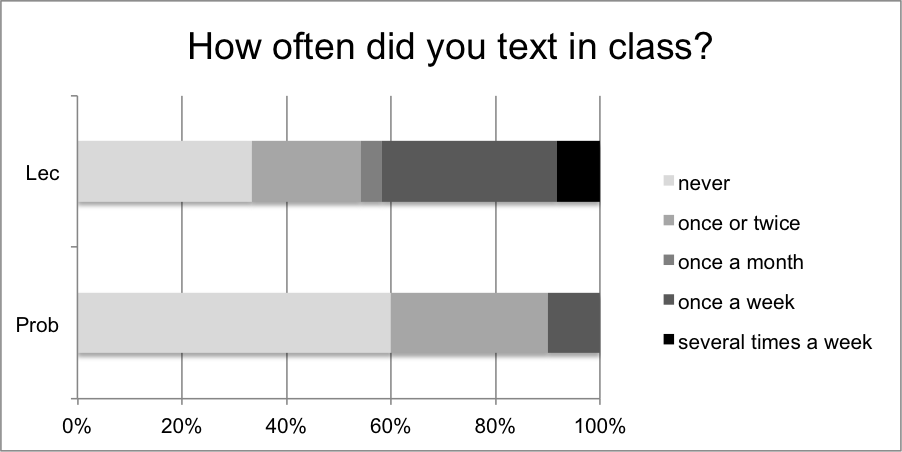}
\end{wrapfigure}
We asked this question at the end of the course to measure intellectual engagement during class, hoping that in a course where students are presenting their solutions to each other and discussing their ideas for all of class time, students would be more actively engaged than in a lecture setting. This was the case; students in the lecture-based sections texted more (difference {\bf is} statistically significant). As the instructor, the only times I saw students using their phones in the problem-solving section was in the down time at the beginning of class when some students had finished writing their solutions on the board and were waiting for the other students to finish.

\subsection*{Other than the mathematical content (triple integrals, Green's Theorem, etc.), what other skills or life lessons did you learn in Math 290-3 this quarter?}
We asked this question at the end of the quarter because college courses should teach more than the syllabus; they should also help students to grow more broadly. For example, Northwestern University's mission statement, which echoes those of other colleges and universities, reads: \textit{Northwestern is committed to excellent teaching, innovative research, and the personal and intellectual growth of its students in a diverse academic community.} 

Almost all students answered this question, reporting at least one non-mathematical skill that they had developed in the course. Students in all three sections reported that they had learned better time management, and had learned to explain math to others. 

Comments from students in the problem-solving section had broader applicability to life than comments from students in the lecture-based sections. Typical responses from the lecture-based sections were ``Start early on homework assignments'' and ``Explaining my work and reasoning in Salons.'' These responses indicate intellectual growth. Typical responses from the problem-solving section were  ``Speaking in front of an audience, being prepared, perseverance'' and ``Problem-solving, discussing solutions and different ways to approach problems, asking relevant questions.'' These responses indicate both intellectual and personal growth, and perhaps appreciation of a diverse academic community, regarding different students' ways of approaching problems.

\clearpage
\subsection*{Which of the two statements below more accurately reflects your thoughts regarding the method in which your Math 290-3 course was taught this quarter?}
\begin{wrapfigure}[13]{R}{0.6\textwidth}
\centering
\vspace{-1.5em}
\includegraphics[width=0.6\textwidth]{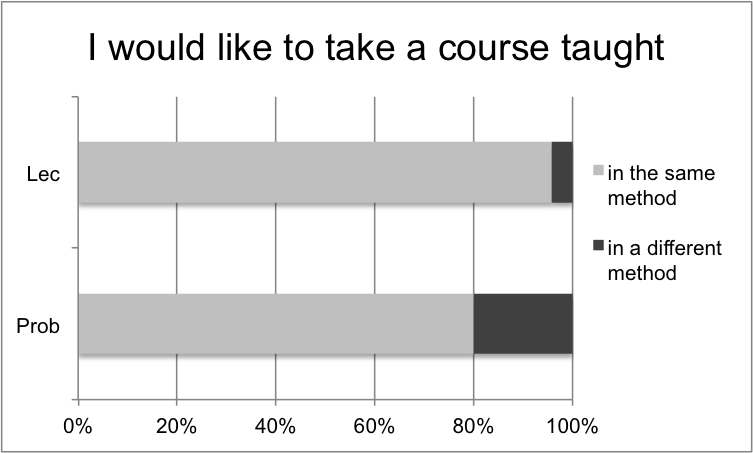}
\end{wrapfigure}
We asked this question because this course was run as a pilot, so the math department would like to know whether it should run more courses like this in the future. Across all sections, a majority of students wished to take another course taught using the same method as the one that they were taking. 

The percentage of students wishing to take a course using a different method was higher in the problem-solving section than in the lecture-based sections (difference is not statistically significant). Based on student evaluations, this was probably because they felt like it took too much time to prepare for class. In my opinion, it took the same amount of time, but spread more evenly across the semester: I expected students to spend about $8$ hours per week on the course, two hours doing the homework for each of 4 class periods. Since the problems constantly required using all the tools developed up to this point, and since the homework and discussion for two class periods before the exam were practice exams, there was no need for additional studying for an exam. In contrast, in the lecture sections (in all three quarters) students only spent about 5 hours per week on the weekly problem set, but then 10 or more hours studying for each exam. 

\section{Long-range outcomes}\label{longrange}

A year after the completion of Math 290-3, we collected long-range data from the registrar. We determined how many math courses each student had taken in their sophomore year of college. We also determined how many of the students had decided to declare the math major. In the table below, we report the number of courses taken by each student, as well as the averages and standard deviations, and the math major declaration data, for each of the three sections. 

\begin{table}[!h]
\centering
\label{scores}
\begin{tabular}{c|cc|l}
 & $\#$ of courses taken & Average &  proportion declaring math major  \\ \hline
Problem-solving & 2, 3, 3, 3, 4, 5, 5, 6, 7 & 4.2  & 7/10 = 70$\%$ \\ \hline
Lecture-based 1  & 0, 1, 1, 1, 1, 3, 3, 4, 5, 6  & 2.5  & 5/10 = 50$\%$  \\  
Lecture-based 2  &  0, 0, 1, 1, 1, 1, 1, 3, 4, 6, 6, 6 &   2.5  & 6/13 = 46$\%$
\end{tabular}
\end{table}

Every student who was in the problem-solving section in their freshman spring took at least two math courses in their sophomore year.\footnote{Only $9$ course counts are reported for the problem-solving section because one of the $10$ students was not enrolled in the academic year 2016-2017. This student did declare the math major.} By contrast, half of the students from each lecture section in the freshman spring took at most 1 math course in their sophomore year. 

Students in the problem-solving section averaged 4.2 math courses in their sophomore year, while students in the lecture sections averaged just 2.5 (difference {\bf is} statistically significant). Additionally, $70\%$ of students in the problem-solving section declared the math major, while about $50\%$ of students in each lecture section did so (difference is not statistically significant). 

While conclusions of causation are confounded by the fact that the problem-solving section was opt-in (see $\S$\ref{structure}), this data does show that students in the problem-solving course chose to engage more with mathematics over the following year than did those in the lecture sections. It is striking that the numbers are almost exactly the same for the two lecture sections, despite their having different instructors. 

\section{Applicability of the method to less ideal situations}\label{applications}

The Math 290-3 course was the best-case scenario for trying a new method: The class was very small, the student population was homogeneous (see $\S\ref{structure}$), and the students voluntarily chose to take a course under the problem-based method. Here I will discuss adaptations to other settings.

\vspace{0.05in}\noindent{\bf 15 students, heterogeneous bridge course} At the same time that I was teaching Math 290-3, I taught the Northwestern math department's Introduction to Mathematical Proof course using the problem-centered, discussion-based method. Students did homework problems, and class time consisted of students presenting and discussing their solutions, just as in Math 290.

Mine was one section of a multi-section course, where the other sections were taught in a different style. Students did not know that my section would be taught in this method when they signed up, so unlike the registration process for Math 290, the registration for this course was {\it not} opt-in. Nonetheless, not a single student even inquired about switching to another section. When I met with the students individually, many of them said that they really liked the way we spent class time and would not want anything changed. The only substantive change any student suggested was to request that I hand out full solutions to all of the assigned problems. 

\begin{figure}[!h]
\noindent\includegraphics[trim={1cm 0 1cm  6cm},clip,width=\textwidth]{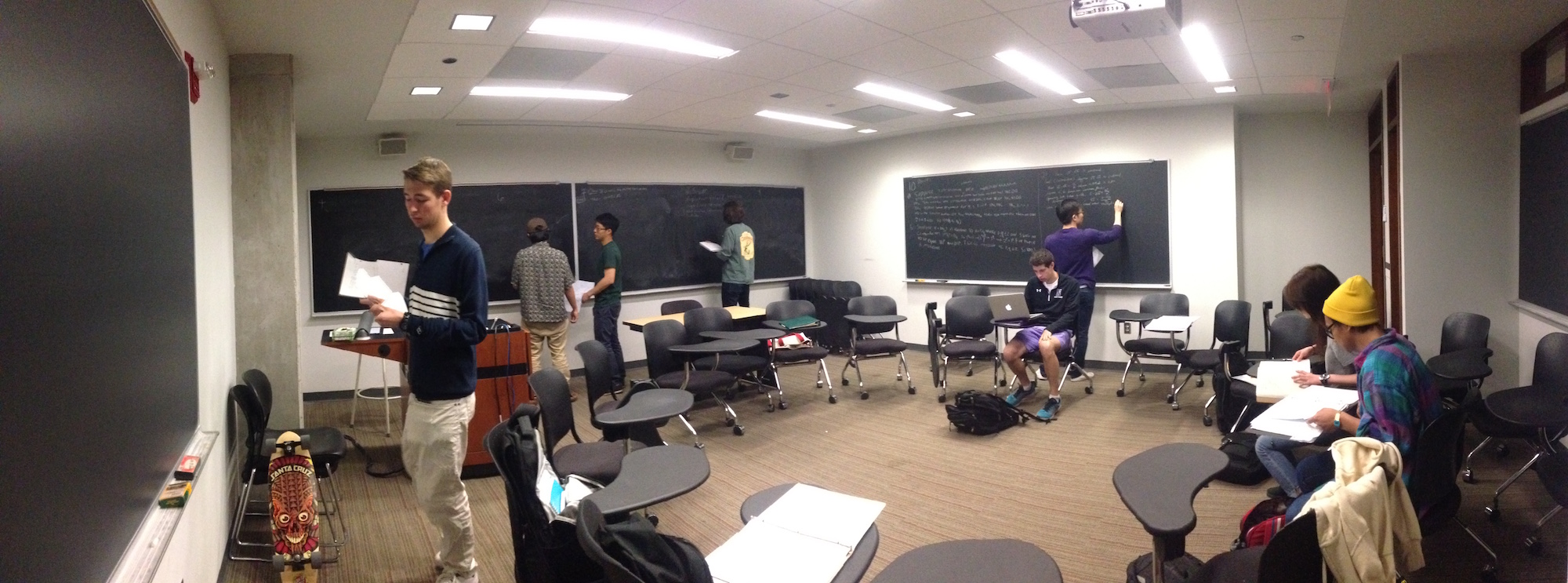}
\captionsetup{labelformat=empty}
\caption{Math 300 students writing up problems in the first $5$ minutes of class, with circle of chairs}
\end{figure}

There were $15$ students in the Math 300 course, and the class was more heterogeneous than the Math 290 course, both of which required modifications to the Harkness method. With the larger class size, each student participated for a smaller fraction of class time, and it was easier for tentative students to hide. The class had students from all four class years, with several of the students winning prizes from the math department for excellence, while others struggled to grasp the material. This made it more difficult to set the pace of the class, because the weaker students needed to spend time discussing the foundational problems, while the stronger students wanted to move on and get to the more challenging problems. 

To ensure that everyone spent time thinking hard about math at an appropriate level, I ended up having the students spend about $30$ minutes presenting and discussing problems as in the Math 290 course, and the remaining $15$ minutes in pairs at the board working on homework problems that we had not discussed, so that all would be forced to think, work and talk.

\vspace{0.05in}\noindent{\bf 25 students, heterogeneous senior seminar} In the fall of 2016, I taught a senior seminar at Williams College in my research area of mathematical billiards. I wrote my own problem-based curriculum \cite{billiards}, and the students did problems for homework and discussed them in class as described above. Some of the students were advanced juniors, some were seniors planning to go to graduate school in math, and some were relatively weak seniors who just needed to take the course to fulfill the ``400-level seminar'' requirement of the math major. When students pre-registered for the course in the spring, it was described as a lecture course, and I sent them an email (including an early version of this paper, which many of them read) in the late summer telling them that it would be taught in this discussion-based, problem-centered format.

In this case, 25 students is far too many for a single discussion, so when they came in they chose a playing card. A black card assigned them to the blackboard group, and red to the whiteboard group. I wrote the problem numbers on each of the two boards, students wrote up full solutions on each board, and each discussion operated independently (see picture).

\begin{figure}[!h]
\noindent\includegraphics[width=\textwidth]{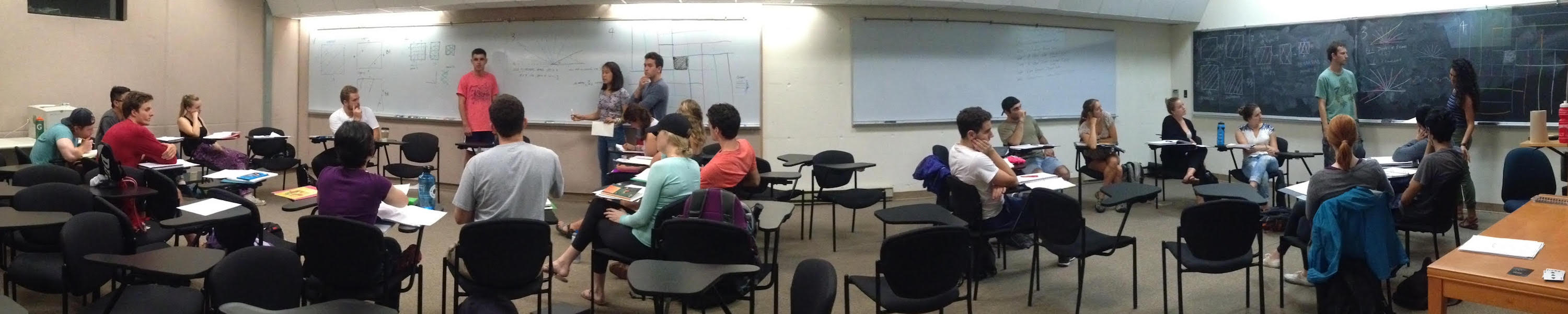}
\captionsetup{labelformat=empty}
\caption{Two simultaneous discussions in the billiards senior seminar at Williams College}
\end{figure}

I sat between the two groups and gave help or clarification when needed. Students were distracted by the noise of the other group for the first few minutes of class, but as soon as they got interested in their discussion, they were focused on their own group. In response to student feedback, we tried 3 groups of 8 students and 4 groups of 6 students, but in that case there was too much noise in the room, and frequently there was not enough critical mass of students in each group to solve the hard problems, or to ask enough questions to create a rich discussion.

Students were very positive about the course format. Some of the senior math majors said that they had never spoken before in a college math course, and now they were speaking every day. They pointed out a few drawbacks of having two simultaneous discussions, such as the noise from the other group, and that sometimes one group needed help when I was busy working with the other group, but on the whole they thought the method was better than lectures.

\vspace{0.05in}\noindent{\bf 13 students, high school classes}  As mentioned in $\S\ref{design}$, Phillips Exeter Academy, a private four-year boarding high school in New Hampshire, developed this type of problem-based curriculum, and the discussion-based ``Harkness method''  \cite{book,harkness}. They use the Harkness method in {\it every} class in all departments, for {\it every} student, and they use a problem-based curriculum in {\it all} math classes. Their experience is that when students first take a problem-based math class, about $70\%$ of students like it and $30\%$ dislike it. After they become acclimated to it, about $90\%$ of students end up liking it, and about $10\%$ of students persist in their dislike of the method. 

\newpage
\section{Conclusions}
On every survey question that we asked the Math 290-3 students, students in the problem-solving section reported higher self ratings, and their responses showed a larger increase in average score, than those of students in the lecture-based sections. Some of these differences are not statistically significant, but this is expected with such small sample (class) sizes, with only $10$ in the problem-based class. The fact that there are some statistically significant (positive) results with such small numbers demonstrates the educational power of inquiry-based learning. Certainly, larger sample sizes and more data would allow more definitive statistical inferences to be drawn.

Additionally, the outcomes of students in the problem-solving section were compared against the outcomes of students in lecture-based sections with excellent instructors who used active learning techniques, as discussed in $\S\ref{instructors}$. Still, the students in the problem-solving section had comparable exam performances,  reported higher self ratings on the other ``soft skills'' measures, and took more math courses the following year, compared to students in these excellent lecture courses.

As explained in $\S$\ref{structure}, all of the students in the problem-solving section opted in, which gave us a ``best-case scenario,'' where all of the students were invested in the method. People who sat in on the course noted that the group had a good rapport, and that these students seemed well-suited to discussing mathematics.  
One might argue that the method would not work as well if the problem-based method were forced on students who had no choice in the matter. However, students were quite happy with the course in each of the three other class situations discussed in $\S\ref{applications}$, which were not offered on an opt-in basis but instead were required for all students.

This type of problem-based approach is neither a teacher-proof curriculum, nor a curriculum-proof method: its success depends on the skills of the instructor, on the particular problems chosen for homework and discussion, and on the students' coming to class inquisitive and prepared. Too much talking from the instructor turns the class silent, while not enough guidance can allow the students to confuse each other. But with a good problem set, and with good preparation of the instructor, the class can be wonderful, as it was here. 

My dream is to lead a revolution, to a world where college math classes are crafted to create opportunities for student discovery, are centered around student inquiry, and are built on authentic human interactions among students and with the instructor. 
I encourage you the reader, and indeed all instructors of mathematics at all levels, to try a problem-solving curriculum, and a discussion-based course, to whatever extent you can, in your math courses.

\vspace{0.5in}

\subsection*{Acknowledgements}
Paul Goerss, Ursula Porod and John Alongi in the mathematics department at Northwestern University approved my plan to run an experimental section of Math 290, allowing me to live out my dream. Emily Kefferstan made the opt-in registration happen smoothly and collected the long-range student data. Aaron Peterson and Martha Precup taught the lecture-based sections of the course, ensured that all of their students took all of the surveys, and were endlessly accommodating and helpful with this endeavor. Paul VanKoughnett was an excellent graduate teaching assistant, solving all of the problems and suggesting many good changes before we gave them to the students. Paul, along with Keith Burns, taught the class for a week while I was at a conference. At Phillips Exeter Academy, Rick Parris wrote the majority of the Exeter mathematics problems that I used for this course, and Gwyn Coogan taught me how to write a problem-based curriculum customized for my students. Laura Batterink analyzed all of the data to determine which results were statistically significant. The two anonymous referees made suggestions that substantially improved this paper. Most of all, the $10$ students who took the problem-solving course bought in fully from the very beginning, believed in what we were doing, and were a joy to work with every single day.


\vfill

{\footnotesize \bibliography{thesis}}

\end{document}